\documentclass[preprint,12pt]{elsarticle}
\usepackage{amssymb,graphicx}
\newtheorem{theorem}{Theorem}[section]
\newtheorem{lemma}[theorem]{Lemma}
\newtheorem{definition}[theorem]{Definition}
\newtheorem{example}[theorem]{Example}
\newtheorem{proposition}[theorem]{Proposition}

\newtheorem{corollary}[theorem]{Corollary}
\newtheorem{remark}[theorem]{Remark}

\journal{Fund. Math. in about December 2012.}

\begin{document}

\begin{frontmatter}

\title{When will the crossing number of an alternating link decrease by two via
a crossing change?}

\author{Xian'an Jin\footnote{Email:xajin@xmu.edu.cn},\ \ Fuji Zhang\footnote{Email:fjzhang@xmu.edu.cn},\ \ Jun Ge}
\address{School of Mathematical Sciences,\\
Xiamen University, Xiamen, Fujian 361005,\\
P. R. China\\
This is the revised version.}

\begin{abstract}
Let $D$ be a reduced alternating diagram of a non-split link $L$ and
$\tilde{L}$ be the link whose diagram is obtained from $D$ by a
crossing change. If $\tilde{L}$ is alternating, then
$c(\tilde{L})\leq c(L)-2$. In this paper we explore when
$c(\tilde{L})=c(L)-2$ holds and obtain a simple sufficient and
necessary condition in terms of plane graphs corresponding to $L$.
This result is obtained via analyzing the behavior of the Tutte polynomial of the signed plane graph corresponding to $\tilde{L}$.
\end{abstract}

\end{frontmatter}

\section{Introduction}

\begin{figure}[htbp]
 \centering
\includegraphics[width=5cm]{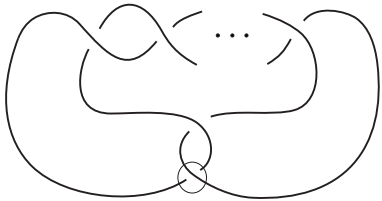}\label{unknotting}
\renewcommand{\figurename}{Fig.}
\caption{By changing the crossing circled, the alternating knot is
unknotted.}
\end{figure}

Let $L$ be a link. We denote by $c(L)$ the \emph{crossing number} of
the link $L$, that is, the smallest number of crossings, the minimum
being taken over all diagrams of $L$. Let $D$ be a diagram, by
\emph{a crossing change} we mean exchanging the over-pass and the
under-pass curves at a single crossing of $D$. The crossing number of an
alternating link may decrease dramatically via a single crossing
change, for example, alternating knots with unknotting number one as
shown in Fig. \ref{unknotting}. It is natural to ask which
conditions should be satisfied by an alternating knot diagram such
that its crossing number decreases only a little when one changes
any its crossing. It is well known that there is a one-to-one correspondence between link
diagrams and signed plane graphs via the medial construction, which
provide a method of studying knots using graphs \cite{Ad}.
We shall answer this question in terms of corresponding plane graphs
under some moderate conditions.

Another inspiration for our study is works of ordering knots via
crossing changes. In \cite{Diao}, Diao et al defined a partial
ordering of links using a property derived from their minimal
diagrams. A link $L'$ is called a predecessor of a link $L$ if
$c(L')<c(L)$ and a diagram of $L'$ can be obtained from a minimal
diagram $D$ of $L$ by a single crossing change. In addition, in
\cite{Ta}, Taniyama defined that $L_1$ is a major of $L_2$ if every
diagram of $L_1$ can be transformed into a diagram of $L_2$ by
applying crossing changes at some crossings of the diagram of $L_1$.
The notion of major is extended to $s$-major in \cite{Endo} via
adding smoothing operations by Endo et al. Our result may help to their studies.

We noted the following result obtained by L. Wu et al.
\begin{theorem}\cite{Lei}\label{lei}
Let $L$ be a non-split link which admits a reduced alternating
diagram $D$. Let $\tilde{L}$ be the link obtained from $D$ by a
crossing change. If $\tilde{L}$ is alternating, then
\begin{eqnarray}
c(\tilde{L})\leq c(L)-2.
\end{eqnarray}
\end{theorem}
Theorem 3.2 in \cite{Diao} shows that Theorem \ref{lei} holds for
rational links. In this paper we shall explore when the equality in
Theorem \ref{lei} holds, that is, when the crossing number of an
alternating link decreases by two via a crossing change?

We attempt to study the effect of crossing number of a link after a
single crossing change and find that it is difficult to deal with it by using the diagrammatic approach.
However, when we turn to the corresponding plane graphs, the Tutte polynomial of graphs or signed graphs
provides a good tool to solve the problem.
Let $G$ be a graph. The multiplicity $\mu(e)$ of an edge $e=(u,v)$ of $G$ is
the number of all edges with end-vertices $u$ and $v$. We use $N(v)$
to denote the set of all vertices of $G$ that have a common edge with $v$. In this
paper we proved

\begin{theorem}\label{eq} Let $G$ be a connected bridgeless and loopless
positive plane graph and $e=(u,v)$ be an edge of $G$. Let $L$ be the
alternating link corresponding to $G$ and $\tilde{L}$ be the link
corresponding to $\tilde{G}$ obtained from $G$ by changing the sign
of $e$ from $+$ to $-$. Suppose that $L$ is non-split and
$\tilde{L}$ is alternating, we have
\begin{enumerate}
\item if $\tilde{L}$ is split, then
$c(\tilde{L})=c(L)-2$ if and only if $\mu(e)=2$ and if we suppose
that $f$ is the edge parallel to $e$, then $G-e-f$ is disconnected.
\item if $\tilde{L}$ is non-split, then $c(\tilde{L})=c(L)-2$ if and
only if one of the following two
conditions holds:
\begin{enumerate}
\item[(1)] $\mu(e)=1$, $G-e$ has bridges and $N(u)\cap N(v)=\emptyset$.
\item[(2)] $\mu(e)>1$ and if we suppose that $f$ is an edge parallel to
$e$, then $G-e-f$ is connected and bridgeless.
\end{enumerate}
\end{enumerate}
\end{theorem}

Note that the characterization of plane graphs corresponding to
knots has been given in \cite{Sh,EP,Go}. In the following of this
section, we apply Theorem \ref{eq} to the case of knots. A graph is
said to be \emph{2-edge connected} if it is connected and
bridgeless. An edge with multiplicity 1 or a (not necessarily
maximal) multiple edge, which is formally defined in Section 4, of a
2-edge connected graph $G$ is said to be \emph{reducible} if $G$ is
still 2-edge connected after deleting the edge or the multiple edge,
otherwise it is said to be irreducible. A triangle in a graph $G$ is
called to be \emph{quasi-simple} if it has at least one edge with
multiplicity 1.

\begin{corollary}\label{eq1}
Let $G$ be a connected bridgeless and loopless positive plane graph.
Let $L$ be the alternating link corresponding to $G$ and $\tilde{L}$
be any link corresponding to $\tilde{G}$ obtained from $G$ by
changing the sign of an edge of $G$ from $+$ to $-$. Suppose that
$L$ is a knot and $\tilde{L}$ is always alternating. Then
$c(\tilde{L})=c(L)-2$ for any $\tilde{L}$ if and only if
\begin{enumerate}
\item $G$ is quasi-simple triangle free,
\item each edge with multiplicity 1 is irreducible,
\item A pair of edges in any maximal multiple edge is reducible.
\end{enumerate}
\end{corollary}

\noindent{\bf Proof.} Since $L$ is a knot, $\tilde{L}$ is also a
knot. Hence both $L$ and $\tilde{L}$ are non-split. Let $e$ be an
edge of $G$. If $\mu(e)=1$, Conditions 1 and 2 are equivalent to
Theorem \ref{eq} 2(1). If $\mu(e)>1$, Condition 3 is equivalent to
Theorem \ref{eq} 2(2). $\Box$

A 2-edge connected graph $G$ is said to be \emph{minimal} if, for
each edge $e$ of $G$, $G-e$ has bridges. We further restrict
ourselves to simple graphs, that is, graphs having no loops or
multiple edges, and, as a direct consequence of Corollary \ref{eq1},
we obtain

\begin{corollary}\label{eq2}
Let $G$ be a connected bridgeless and loopless positive simple plane
graph. Let $L$ be the alternating link corresponding to $G$ and
$\tilde{L}$ be any link corresponding to $\tilde{G}$ obtained from
$G$ by changing the sign of an edge of $G$ from $+$ to $-$. Suppose
that $L$ is a knot and $\tilde{L}$ is always alternating. Then
$c(\tilde{L})=c(L)-2$ for any $\tilde{L}$ if and only if $G$ is a
triangle-free and minimal 2-edge connected graph.
\end{corollary}

Compared with Theorem \ref{eq}, Corollaries \ref{eq1} and \ref{eq2}
can both be viewed as results on the 'whole' alternating link
diagram. For the construction and properties of minimal 2-edge
connected graph, see \cite{Zhu,Chache}.

The paper is organized as follows. In Section 2, we provide some
preliminary knowledge, including the relation between the crossing
number of an alternating link and the span of its Jones polynomial,
and the relation between the Jones polynomial and the Tutte
polynomial. We then give a graph-theoretic proof of Theorem
\ref{lei} in Section 3. In Section 4, we
obtain a 'dual' result of Dasbach and Lin \cite{Lin} on the
coefficients of $T_G(-t,-t^{-1})$. Theorem \ref{eq} is thus obtained
by studying the proof in Section 3 and using the 'dual' result and
its proof is given in Section 5. In the final Section 6, we give an
example illustrating Theorem \ref{eq} and pose two problems for further study.

\section{Preliminaries}

The readers who are familiar with the knowledge on the correspondence
between graphs and links, Jones polynomial and Tutte polynomial can skip this section.

\subsection{Some terminologies and notations}

A \emph{graph} $G$ is a pair of sets $V(G)$ and $E(G)$, where $V(G)$
is a non-empty finite set (of \emph{vertices}) and $E(G)$ is a
multi-set of unordered pairs $(x,y)$ (not necessarily distinct) of
vertices called \emph{edges}. An edge with unordered pair $(x,x)$ is
called a \emph{loop}. For $v\in V(G)$, let $N(v)=\{u\in V(G)|(u,v)\in E(G)\}-\{v\}$.
Graphs can be represented graphically,
that is, we can draw it as follows: each vertex is indicated by a
point, and each edge $(x,y)$ by a line joining the points $x$ and
$y$. A graph is \emph{planar} if it can be embedded in the plane,
that is, it can be drawn on the plane so that no two edges
intersect. A \emph{plane} graph is a particular plane embedding of a
planar graph. A graph is said to be \emph{trivial} if it consists of
only an isolated vertex without loops. A \emph{signed} graph is a graph each of
whose edges is labeled with a sign ($+$ or $-$).

A graph is said to be \emph{connected} if, for any its two distinct
vertices $u,v$, there is a path $u=u_0u_1u_2\cdots u_l=v$, where
$u_i$ ($i=0,1,\cdots,l$) are all distinct and $(u_{i-1},u_i)$ is an
edge for $i=1,2,\cdots,l$. A \emph{connected component} of a graph is a maximal
 connected subgraph of the graph. A \emph{bridge} of a graph $G$ is an edge
whose removal would increase the number of connected components of
$G$. By contracting an edge we mean deleting the edge firstly and
then identifying its end-vertices. Let $e$ be an edge of $G$. We
shall denote by $G-e$ and $G/e$ the graph obtained from $G$ by
deleting and contracting the edge $e$, respectively. When $G$ is a
plane graph, $G-e$ and $G/e$ are also plane graphs obtained in a
natural way.

A \emph{knot} is a simple closed piecewise linear curve in Euclidean
3-space $\mathbb{R}^3$. A \emph{link} is the disjoint union of
finite number of knots, each knot is called a \emph{component} of
the link. We take the convention that a knot is a one-component
link. We can always represent links in $R^3$ by \emph{link diagrams}
in a plane, that is, regular projections with a short segment of the
underpass curve cut at each double point of the projection.

A link diagram is said to be \emph{split} if it is a composition of
the diagrams of two links with no points in common \cite{M}, and
otherwise \emph{non-split} or \emph{connected}. A link that has a
split diagram is said to be a \emph{split} link, and otherwise
\emph{non-split} or \emph{connected}. A link diagram is said to be
\emph{alternating} if over- and under-crossings alternate as one
travels the link (crossing at the crossings), and otherwise
\emph{non-alternating}. A link is said to be \emph{alternating} if
it has an alternating link diagram, and otherwise
\emph{non-alternating}. A \emph{nugatory} crossing of a link diagram
is a crossing in the diagram so that two of the four local regions
at the crossing are part of the same region in the larger diagram. A
\emph{reduced} diagram is one that does not contain nugatory
crossings.

\subsection{Links and graphs}

The 1-1 correspondence between link diagrams and signed plane graphs
has been known for about one hundred years. It was once one of the
methods used by Tait and Little in the late 19th century to
construct a table of knot diagrams of all knots starting with graphs
with a relatively small number of edges and then increasing the
number of edges \cite{M}. To describe this correspondence, we first
recall the medial graph of a plane graph.

\begin{definition}
The medial graph $M(G)$ of a non-trivial connected plane graph
$G$ is a 4-regular plane graph obtained by inserting a vertex on
every edge of $G$, and joining two new vertices by an edge lying in
a face of $G$ if the vertices are on adjacent edges of the face; if
$G$ is trivial, its medial graph is defined to be a simple closed
curve surrounding the vertex (strictly, it is not a graph); if a
plane graph $G$ is not connected, its medial graph $M(G)$ is defined
to be the disjoint union of the medial graphs of all its connected
components.
\end{definition}

\begin{figure}[htbp]
\centering
\includegraphics[width=7cm]{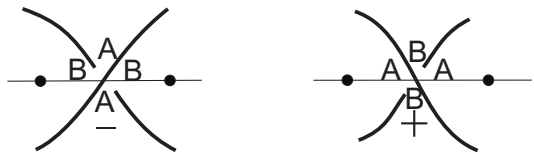}
\renewcommand{\figurename}{Fig.}
\caption{The correspondence between a crossing and a signed
edge \cite{K3}.}\label{cors}
\end{figure}

Given a signed plane graph $G$, we first draw its medial graph
$M(G)$. To turn $M(G)$ into a link diagram $D(G)$, we turn the
vertices of $M(G)$ into crossings by defining a crossing to be over
or under according to the sign of the edge as shown in Fig.
\ref{cors}. Conversely, given a connected link diagram $D$, shade it
as in a checkerboard so that the unbounded face is unshaded. Note
that such a shading is always possible, since link diagrams can be
viewed as 4-regular plane graphs, see Exercise 9.6.1 of
\cite{Bondy}. We then associate $D$ with a signed plane graph $G(D)$
as follows: For each shaded face $F$, take a vertex $v_F$, and for
each crossing at which $F_1$ and $F_2$ meet, take an edge
$(v_{F_1},v_{F_2})$ and give the edge a sign also as shown in Fig.
\ref{cors}. if a link diagram $D$ is not connected, its
corresponding signed plane graph $G(D)$ is defined to be the
disjoint union of the signed plane graphs of all its connected
components.

Under the 1-1 correspondence described above, there is also an 1-1
correspondence between crossings of $D$ and edges of $G(D)$. The
following three properties on the correspondence are all obvious.

\noindent {\bf P1:}\ $D$ is a connected link diagram if and only if
its corresponding signed plane graph $G(D)$ is connected.

\noindent {\bf P2:}\ A crossing of $D$ is nugatory if and only if
its corresponding edge in $G(D)$ is a loop or a bridge. Furthermore,
$D$ is reduced if and only if $G(D)$ is loopless and bridgeless.

\noindent {\bf P3:}\ $D$ is alternating if and only if all edges of
$G(D)$ have the same signs.

\subsection{Jones and Tutte polynomials}

Let $L$ be an oriented link, $V_L(t)$ be the Jones polynomial
\cite{Jones} of $L$. We denote by $span_v(L)$ the difference between
the maximal and minimal degrees of $V_L(t)$, i.e.
\begin{eqnarray*}
span_v(L)=\max \deg V_L(t)-\min \deg V_L(t).
\end{eqnarray*}

In \cite{K1}, Kauffman introduced the Kauffman bracket polynomial of
unoriented link diagrams. Let $D$ be an unoriented link diagram. Let
$[D]=[D](A,B,d)$ be the Kauffman square bracket polynomial of $D$,
$<D>\ =[D](A,A^{-1},-A^2-A^{-2})$ be the Kauffman bracket polynomial
of $D$. We denote by $span_k(D)$ the difference between the maximal
and minimal degrees of $<D>$, i.e.
\begin{eqnarray*}
span_k(D)=\max \deg <D>-\min \deg <D>.
\end{eqnarray*}

Let $L$ be an oriented link, $D$ be an oriented diagram of $L$. The
writhe $w(D)$ of $D$ is defined to be the sum of signs of the
crossings of $D$. Kauffman proved \cite{K1,K2}
\begin{eqnarray*}
V_L(t)=(-A^3)^{-w(D)}<D>|_{A=t^{-1/4}}.
\end{eqnarray*}
Hence we have
\begin{eqnarray}
span_v(L)={1\over 4}span_k(D).
\end{eqnarray}

\begin{lemma}\cite{K1,Mu,Th}\label{span} Let $D$ be a unoriented diagram of an oriented link
$L$.
\begin{enumerate}
\item[(1)] If $L$ is a non-split alternating link, then $c(L)=span_v(L)={1\over 4}span_k(D)$.
\item[(2)] If $L$ is a split alternating link with $n(L)$ non-split components, then
$c(L)=span_v(L)-n(L)+1={1\over 4}span_k(D)-n(L)+1$.
\end{enumerate}
\end{lemma}

Given a crossing of a link diagram, we can distinguish two out of the four small regions incident at the crossing. Rotate
the over-crossing arc counterclockwise until the under-crossing arc is reached, and call the small two regions swept out the \emph{$A$-channels} and other two the \emph{$B$-channels}.
For example, in Fig. \ref{cors}, the edge with sign $+$(resp. $-$) edge crosses $A$-channels (resp. $B$-channels). In the case of an alternating link diagram, each of its regions has only $A$-channels or only
$B$-channels. Calling a region an \emph{$A$-region} if all its channels are $A$ channels, and a \emph{$B$-regio}n if all its channels are $B$ channels.

\begin{lemma}\cite{K1,K2}\label{degree}
Let $D$ be a connected reduced alternating link diagram. Then
\begin{enumerate}
\item[(1)] $\max \deg <D>\ =V+2W-2$, where $V$ is the number of crossings of
$D$ and $W$ is the number of $B$-regions. The coefficient of this
power of $A$ in $<D>$ is $(-1)^{W-1}$.
\item[(2)] $\min \deg <D>\ =-V-2B+2$, where $V$ is the number of crossings of
$D$ and $B$ is the number of $A$-regions. The coefficient of this
power of $A$ in $<D>$ is $(-1)^{B-1}$.
\end{enumerate}
\end{lemma}

Motivated by the 1-1 correspondence between link diagrams and signed
plane graphs, in \cite{K3} Kauffman constructed a Tutte polynomial
for signed graphs, which is generalizations of both the Tutte polynomial
\cite{Tutte} for ordinary graphs and the Kauffman square bracket polynomial.
Let $G$ be a signed graph and
$Q[G]=Q[G](A,B,d)$ be the Tutte polynomial of $G$, which we shall
call the $Q$-polynomial for clarity.

\begin{definition}\label{d2}
The $Q$-polynomial can be defined by the following recursive rules:
\begin{enumerate}
\item Let $E_n$ be the edgeless graph with $n$ vertices. Then
    \begin{eqnarray*}
     Q[E_n]=d^{n-1}.
    \end{eqnarray*}
\item
   \begin{enumerate}
     \item If $e$ is a bridge, then
      \begin{eqnarray*}
            & &Q[G]=(A+Bd)Q[G/e] \  { when} \  s(e)=+ \ and\\
            & &Q[G]=(B+Ad)Q[G/e] \  {when} \  s(e)=-.
      \end{eqnarray*}
    \item If $e$ is a loop, then
      \begin{eqnarray*}
            & &Q[G]=(B+Ad)Q[G-e] \ { when} \  s(e)=+ \ and\\
            & &Q[G]=(A+Bd)Q[G-e] \ { when} \  s(e)=-.
      \end{eqnarray*}
     \item If $e$ is neither a bridge nor a loop, then
      \begin{eqnarray*}
      & &Q[G]=BQ[G-e]+AQ[G/e]\ { when} \ s(e)=+ \ and\\
      & &Q[G]=AQ[G-e]+BQ[G/e]\ { when} \ s(e)=-.
       \end{eqnarray*}
    \end{enumerate}
\end{enumerate}
\end{definition}

\begin{lemma}\label{tr}\cite{K3}
Let $G$ be a signed plane graph, $D(G)$ be the link diagram
corresponding to $G$. Then $Q[G]=[D(G)]$.
\end{lemma}

Let $G$ be a signed plane graph. The componentwise dual $G^*$ of $G$
is defined to be the disjoint union of the dual graphs of all
connected components of $G$. Note that there is a bijection between
edges of $G$ and edges of $G^*$, and the edge $e\in E(G)$ and the
corresponding edge $e^*\in E(G^*)$ receive opposite signs.

\begin{lemma}\label{dual}
Let $G$ be a signed plane graph, $G^*$ be the componentwise dual of
$G$. Then $Q[G]=Q[G^*]$.
\end{lemma}

From now on we always suppose that
$Q[G]=Q[G](A,A^{-1},-A^2-A^{-2})$. Recall that the Tutte polynomial
$T_G(x,y)$ of a graph $G=(V,E)$ can be defined by the following
summation:
\begin{eqnarray}
\sum_{F\subseteq E}(x-1)^{k(F)-1}(y-1)^{|F|-|V|+k(F)},
\end{eqnarray}
where $k(F)$ is the number of connected components of the spanning
subgraph $(V,F)$ of $G$.

A signed graph $G$ is said to be positive (resp. negative) if any of
its edges receives a positive (resp. negative) sign.  Using the
Receipe Theorem of the Tutte polynomial \cite{B} or Thistlethwaite
Theorem \cite{Th}, we can deduce

\begin{lemma}\label{this}
Let $G=(V,E)$ be a connected graph, $G_+$ be the positive graph
whose underlying graph is $G$. Then
\begin{eqnarray*}
Q[G_+]=A^{-|E|+2|V|-2}T_G(-A^{-4},-A^4).
\end{eqnarray*}
\end{lemma}

\section{The proof of Theorem \ref{lei}}
Let $L$ be a non-split link which admits a reduced alternating
diagram $D$. Since $L$ is non-split, $D$ must be also connected. Let
$G=G(D)$ be the signed plane graph corresponding to $D$. Without
loss of generality we assume that $G$ is positive. Otherwise, by
Lemma \ref{dual} we shall work on $G^*$. Since $D$ is reduced, $G$
is loopless and bridgeless.

Let $\tilde{L}$ be an alternating link whose diagram $\tilde{D}$ is
obtained from $D$ by a crossing $c$ change. Since $D$ is connected,
$\tilde{D}$ is also connected. Let $\tilde{G}=G(\tilde{D})$ be the
signed plane graph corresponding to $\tilde{D}$.  Then $\tilde{G}$
can be obtained from $G$ by changing the sign of an edge $e$
corresponding to $c$ from $+$ to $-$.

Let $span_q(\tilde{G})=\max \deg Q[\tilde{G}]-\min \deg
Q[\tilde{G}]$. By Lemmas \ref{span} and \ref{tr}, we have
\begin{eqnarray}\label{e1}
c(\tilde{L})&\leq& span_v(\tilde{L})\nonumber \\
&=&{1\over
4}span_q(\tilde{G}).
\end{eqnarray}

By Definition \ref{d2} and note that the sign of the edge $e$ in $G$
(resp. $\tilde{G}$) is positive (resp. negative), we have
\begin{eqnarray*}
& &Q[\tilde{G}]=AQ[G']+A^{-1}Q[G''],\\
& &Q[G]=A^{-1}Q[G']+AQ[G''],
\end{eqnarray*}
where $G'=G-e$ and $G''=G/e$. Hence we obtain

\begin{eqnarray}\label{case1}
Q[\tilde{G}]=A^2Q[G]+(A^{-1}-A^3)Q[G'']
\end{eqnarray}
or
\begin{eqnarray}\label{case2}
Q[\tilde{G}]=A^{-2}Q[G]+(A^{1}-A^{-3})Q[G'].
\end{eqnarray}
Since $G$ is loopless and bridgeless, it is clear that $G'=G-e$ is
loopless and $G''=G/e$ is bridgeless.

\noindent{\bf Case 1.} $G''$ is loopless.

In this case $G''$ is connected, loopless, bridgeless and positive,
hence the link diagram corresponding to $G''$ is connected, reduced
and alternating. Let $H$ be a connected plane graph, we shall use
$v(H), e(H)$ and $f(H)$ to denote the number of vertices, edges and
faces of $H$, respectively. By Lemma \ref{degree}, we have:
\begin{enumerate}
\item $\max \deg Q[G]=\max \deg <D>\ =V+2W-2=e(G)+2f(G)-2$ and the
corresponding coefficient of this power is $(-1)^{f(G)-1}$.

\item $\min \deg Q[G]=\min \deg <D>\ =-V-2B+2=-e(G)-2v(G)+2$ and the
corresponding coefficient of this power is $(-1)^{v(G)-1}$.

\item $\max \deg Q[G'']=e(G'')+2f(G'')-2=e(G)+2f(G)-3$ and the
corresponding coefficient of this power is
$(-1)^{f(G'')-1}=(-1)^{f(G)-1}$.

\item $\min \deg Q[G'']=-e(G'')-2v(G'')+2=-e(G)-2v(G)+5$ and the
corresponding coefficient of this power is
$(-1)^{v(G'')-1}=(-1)^{v(G)}$.
\end{enumerate}
Hence,
\begin{enumerate}
\item $\max \deg A^2Q[G]=e(G)+2f(G)$ and the
corresponding coefficient of this power is $(-1)^{f(G)-1}$.

\item $\min \deg A^2Q[G]=-e(G)-2v(G)+4$ and the
corresponding coefficient of this power is $(-1)^{v(G)-1}$.

\item $\max \deg (A^{-1}-A^3)Q[G'']=e(G)+2f(G)$ and the
corresponding coefficient of this power is $(-1)^{f(G)}$.

\item $\min \deg (A^{-1}-A^3)Q[G'']=-e(G)-2v(G)+4$ and the
corresponding coefficient of this power is $(-1)^{v(G)}$.
\end{enumerate}

Note that the maximal (resp. minimal) degree terms of $A^{2}Q[G]$
and $(A^{-1}-A^{3})Q[G'']$ cancel each other. Therefore, by Eq.
(\ref{case1}), we have

\begin{eqnarray*}
\max \deg Q[\tilde{G}]&\leq& e(G)+2f(G)-4,\\
\min \deg Q[\tilde{G}]&\geq& -e(G)-2v(G)+8.
\end{eqnarray*}
So,
\begin{eqnarray*}
span_q(G)=\max \deg Q[G]-\min \deg Q[G]=2e(G)+2f(G)+2v(G)-4,
\end{eqnarray*}
and
\begin{eqnarray*}
span_q(\tilde{G})&=&\max \deg Q[\tilde{G}]-\min \deg
Q<\tilde{G}>\\
&\leq& 2e(G)+2f(G)+2v(G)-12\\
&=&span_q(G)-8.
\end{eqnarray*}
Hence,
\begin{eqnarray*}
c(\tilde{L})&\leq&{1\over 4}span_q(\tilde{G})\\
&\leq&{1\over 4}span_q(G)-2\\
&=&c(L)-2.
\end{eqnarray*}

\noindent{\bf Case 2.}  $G''$ has loops.

\begin{figure}[htbp]
\centering
\includegraphics[width=8.5cm]{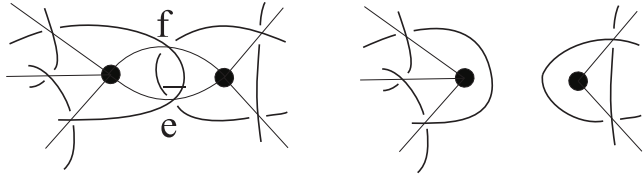}
\renewcommand{\figurename}{Fig.}
\caption{The case $G-e-f$ is disconnected.}\label{2ab}
\end{figure}

Let $f$ be any loop of $G''$. Since $G$ is loopless, $f$ must be an
edge of $G$ parallel to $e$. There are two subcases:

\noindent{\bf Case 2a.} If $G-e-f=G'-f$ is disconnected, then
$\tilde{G}-e-f$ is disconnected. So $\tilde{D}$ can be split as
shown in Fig. \ref{2ab}, which reduces the crossing number by two.
Hence, Theorem \ref{lei} holds.

\noindent{\bf Case 2b.} If $G-e-f=G'-f$ is connected.

Now we prove $G'$ is bridgeless. Firstly $f$ is not a bridge of $G'$
and let $g\neq f$ be an edge of $G'=G-e$. Since $G$ is bridgeless,
$g$ belongs to a cycle $C$ of $G$. If $e\notin E(C)$, $g$ belongs to
a cycle $C$ of $G'$; If $e\in E(C)$, $g$ belongs to a cycle
$C'=C-e+f$ of $G'$. Thus $g$ is not a bridge. Hence $G'$ is
connected, loopless, bridgeless and positive.
Similarly, by Lemma \ref{degree}, we have:
\begin{enumerate}
\item $\max \deg Q[G]=\max \deg <D>\ =V+2W-2=e(G)+2f(G)-2$ and the
corresponding coefficient of this power is $(-1)^{f(G)-1}$.

\item $\min \deg Q[G]=\min \deg <D>\ =-V-2B+2=-e(G)-2v(G)+2$ and the
corresponding coefficient of this power is $(-1)^{v(G)-1}$.

\item $\max \deg Q[G']=e(G')+2f(G')-2=e(G)+2f(G)-5$ and the
corresponding coefficient of this power is
$(-1)^{f(G')-1}=(-1)^{f(G)}$.

\item $\min \deg Q[G']=-e(G')-2v(G')+2=-e(G)-2v(G)+3$ and the
corresponding coefficient of this power is
$(-1)^{v(G')-1}=(-1)^{v(G)-1}$.
\end{enumerate}
Hence,
\begin{enumerate}
\item $\max \deg A^{-2}Q[G]=e(G)+2f(G)-4$ and the
corresponding coefficient of this power is $(-1)^{f(G)-1}$.

\item $\min \deg A^{-2}Q[G]=-e(G)-2v(G)$ and the
corresponding coefficient of this power is $(-1)^{v(G)-1}$.

\item $\max \deg (A^{1}-A^{-3})Q[G']=e(G)+2f(G)-4$ and the
corresponding coefficient of this power is $(-1)^{f(G)}$.

\item $\min \deg (A^{1}-A^{-3})Q[G']=-e(G)-2v(G)$ and the
corresponding coefficient of this power is $(-1)^{v(G)}$.
\end{enumerate}

Note that the maximal (resp. minimal) degree terms of $A^{-2}Q[G]$
and $(A^{1}-A^{-3})Q[G']$ cancel each other. Therefore, by Eq.
(\ref{case2}), we have

\begin{eqnarray*}
\max \deg Q[\tilde{G}]&\leq& e(G)+2f(G)-8,\\
\min \deg Q[\tilde{G}]&\geq& -e(G)-2v(G)+4.
\end{eqnarray*}
Thus,
\begin{eqnarray*}
span_q(\tilde{G})&=&\max \deg Q[\tilde{G}]-\min \deg
Q[\tilde{G}]\\
&\leq& 2e(G)+2f(G)+2v(G)-12\\
&=&span_q(G)-8.
\end{eqnarray*}
Hence,
\begin{eqnarray*}
c(\tilde{L})&\leq&{1\over 4}span_q(\tilde{G})\\
&\leq&{1\over 4}span_q(G)-2\\
&=&c(L)-2.
\end{eqnarray*}
This completes the proof of Theorem \ref{lei}. $\Box$

\section{A 'dual' result}

Let $G=(V,E)$ be a connected loopless graph. $I\subset E$ is said to
be a \emph{multiple edge} if $|I|\geq 2$ and any two of $I$ have the
same end-vertices. A multiple edge $I_M$ is said to be
\emph{maximal} if no multiple edge contains it as a proper subset. In \cite{Lin}, Dasbach abd Lin proved the following lemma.

\begin{lemma}\label{hi}
Let $G=(V,E)$ be a connected loopless graph. Let the Tutte
polynomial evaluation
$$
T_G(-t,-t^{-1})=a_nt^n+a_{n+1}t^{n+1}+\cdots+a_{m-1}t^{m-1}+a_mt^m
$$
with $a_n\neq 0, a_m\neq0$ and $n\leq m$. Then $m=|V|-1$ and

\begin{enumerate}
\item[(1)] $a_m=(-1)^{|V|-1}$.
\item[(2)] $a_{m-1}=(-1)^{|V|-1}(|V|-1-|E|+\sum_{I_M}(|I_M|-1))$, where $I_M$
is a maximal multiple edge and the summation is over all maximal
multiple edges.
\end{enumerate}
\end{lemma}

Let $E_s$ be the edge set of $G_s$, the graph obtained from $G$ by
replacing each maximal multiple edge by a single edge. Then
$a_{m-1}=(-1)^{|V|-1}(|V|-1-|E_s|)$. In the following of this
section, we investigate the value of $n$ and the two coefficients
$a_n$ and $a_{n+1}$, try to obtain a 'dual' result of Lemma
\ref{hi}.

Let $G=(V,E)$ be a connected bridgeless graph. $S\subset E$ is said
to be a \emph{pairwise-disconnecting set} if $|S|\geq 2$ and any two
of $S$ disconnect the graph when deleted. The notion of
pairwise-disconnecting set was introduced in \cite{A}. The following
three statements on pairwise-disconnecting sets are all obvious.

\noindent {\bf ST1:} Any $k$-edge connected graph ($k\geq 3$) does
not contain any pairwise-disconnecting set.

\noindent {\bf ST2:} when $|S|=2$, $S$ is a pairwise-disconnecting
set if and only if $S$ is a 2-edge cut of $G$.

\noindent {\bf ST3:} Any subset with cardinality greater than 1 of a
pairwise-disconnecting set $S$ is also a pairwise-disconnecting set.

\begin{figure}[htbp]
\centering
\includegraphics[width=6.5cm]{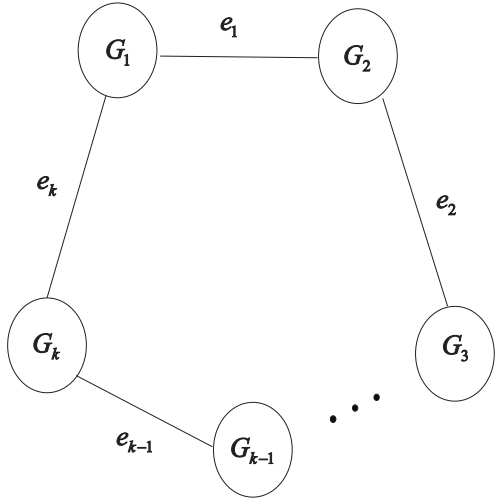}
\renewcommand{\figurename}{Fig.}
\caption{$S=\{e_1,e_2,\cdots,e_k\}$ and $G-S=G_1\cup
G_2\cup\cdots\cup G_k$ and each $G_i$ ($i=1,2,\cdots,k$) is
connected.}\label{cyc}
\end{figure}

\begin{proposition}\label{111}
Let $G=(V,E)$ be a connected bridgeless graph, $S\subset E$ and
$|S|\geq 2$. Then the following are equivalent:
\begin{itemize}
\item $S$ is pairwise-disconnecting set.
\item All edges of $S$ occur on a cycle of $G$ as
shown in Fig. \ref{cyc}.
\item $k(G-S)=|S|$.
\end{itemize}
\end{proposition}

\noindent{\bf Proof.} We first prove that if $k(G-S)=|S|$, then all
edges of $S$ occur on a cycle of $G$ as shown in Fig. \ref{cyc}. It
holds when $|S|=2$ and now we suppose $|S|\geq 3$ and $f\in S$. By
$k(G-S)=|S|$ we have $k(G-S+f)=|S-f|$ and $f$ is a bridge of
$G-S+f$. By induction hypothesis we have $S-f$ occur on a cycle of
$G$. Suppose that $G$ becomes $G_1,G_2,\cdots,G_{|S-f|}$ when $S-f$
deleted. Then $f$ belongs to some $G_i$ and is also a bridge of
$G_i$. Hence, all edges of $S$ occur on a cycle of $G$.

It is clear that if all edges of $S$ occur on a cycle of $G$ as shown
in Fig. \ref{cyc}, then $S$ is a pairwise-disconnecting set.

Finally we prove that if $S$ is pairwise-disconnecting set, then
$k(G-S)=|S|$. It holds when $|S|=2$ and now we suppose $|S|\geq 3$
and $f\in S$. Then $S-f$ is also a pairwise-disconnecting set. By
induction hypothesis we have $k(G-S+f)=|S-f|$. Let $g\in S-f$. Then
$\{f,g\}$ is a 2-edge cut of $G$. Hence, $f$ is a bridge of $G-g$
and also a bridge of $G-S+f$. Therefore, we have
$k(G-S)=k(G-S+f)+1=|S-f|+1=|S|$. $\Box$

\vskip0.2cm

A pairwise-disconnecting set $S_M$ is said to be \emph{maximal} if
no pairwise-disconnecting set contains it as a proper subset.

\begin{proposition}\label{222}
Let $G$ be a connected bridgeless graph. For any given
pairwise-disconnecting set $S$ of $G$, there exists a unique maximal
pairwise-disconnecting set $S_M$ of $G$ containing $S$.
\end{proposition}

\begin{figure}[htbp]
\centering
\includegraphics[width=9.5cm]{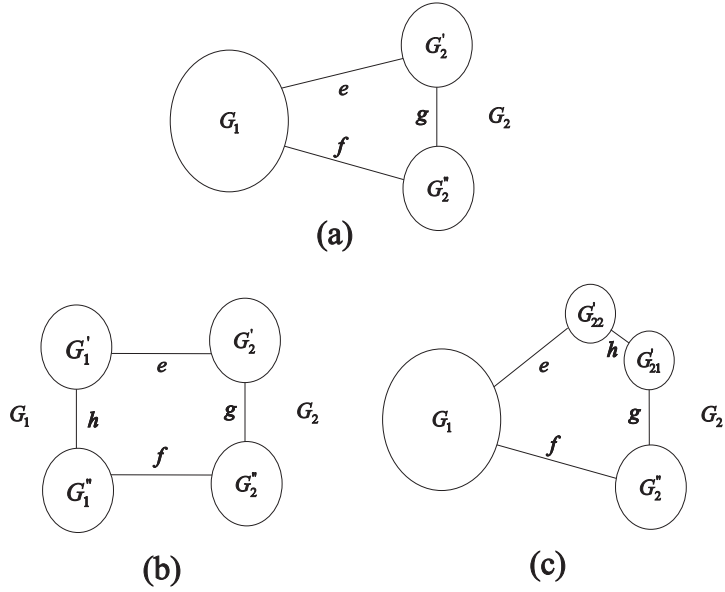}
\renewcommand{\figurename}{Fig.}
\caption{The proof of Proposition \ref{222}.}\label{uni}
\end{figure}

\noindent{Proof.} The existence follows from the definition of
pairwise-disconnecting sets directly. To prove the uniqueness, we
suppose that there are two distinct maximal pairwise-disconnecting
sets $S_M^1$ and $S_M^2$ of $G$ such that $S\subset S_M^i$
($i=1,2$). Let $e,f\in S$ and $g\in S_M^2-S_M^1$. Then $\{e,f\}$ is
a 2-edge cut of $G$ and suppose that $G-e-f=G_1\cup G_2$. Without
loss of generality we suppose that $g\in E(G_2)$. Since $e,f,g\in
S_M^2$ we obtain that $\{e,g\}$ is a 2-edge cut of $G$, which
implies that $g$ must be a bridge of $G_2$. We suppose that
$G_2-g=G_2'\cup G_2''$. See Fig. \ref{uni} (a). Let $h\in
S_M^1-e-f$. We shall show that $\{h,g\}$ is a 2-edge cut of $G$.
Since $e,f,h\in S_M^1$ we obtain that $\{e,h\}$ is a 2-edge cut of
$G$, which implies that $h$ must be a bridge of $G_1$ or $G_2$.
There are two cases.

\noindent {\bf Case 1.} $h$ is a bridge of $G_1$. Suppose that
$G_1-h=G_1'\cup G_1''$, then $G-g-h$ is disconnected as shown in
Fig. \ref{uni} (b).

\noindent {\bf Case 2.} $h$ is a bridge of $G_2$. Without loss of
generality, we suppose that $h\in G_2'$, then $h$ is also a bridge
of $G_2'$. Suppose that $G_2'-h=G_{21}'\cup G_{22}'$, then $G-g-h$
is also disconnected as shown in Fig. \ref{uni} (c).

\noindent Hence, $S_M^1\cup \{g\}$ is a pairwise-disconnecting set,
which contradicts the maximality of $S_M^1$. $\Box$

Now we are in a position to prove a 'dual' result of Lemma \ref{hi}.

\begin{lemma}\label{lo}
Let $G=(V,E)$ be a connected bridgeless graph. Let the Tutte
polynomial evaluation
$$
T_G(-t,-t^{-1})=a_nt^n+a_{n+1}t^{n+1}+\cdots+a_{m-1}t^{m-1}+a_mt^m
$$
with $a_n\neq 0, a_m\neq0$ and $n\leq m$. Then $n=-|E|+|V|-1$ and

\begin{enumerate}
\item[(1)] $a_n=(-1)^{|E|-|V|+1}$.
\item[(2)] $a_{n+1}=(-1)^{|E|-|V|+1}(-|V|+1+\sum_{S_M}(|S_M|-1))$, where $S_M$
is a maximal pairwise-disconnecting set and the summation is over
all maximal pairwise-disconnecting sets.
\end{enumerate}
\end{lemma}

\noindent{\bf Proof.} Recall that
\begin{eqnarray*}
T_G(-t,-t^{-1})&=&\sum_{F\subseteq
E}(-t-1)^{k(F)-1}(-t^{-1}-1)^{|F|-|V|+k(F)}\\
&=&\sum_{F\subseteq
E}(-1)^{|F|-|V|+1}(1+t)^{k(F)-1}(t^{-1}+1)^{|F|-|V|+k(F)}.
\end{eqnarray*}
It is clear that $k(F)-1\geq 0$. Thus we obtain
\begin{eqnarray*}
(1+t)^{k(F)-1}=1+(k(F)-1)t+\left(
\begin{array}{c}
k(F)-1\\
2
\end{array}
\right)t^2+\cdots.
\end{eqnarray*}
Since $|F|-|V|+k(F)$ is the nullity of the subgraph $(V,F)$ of
$G=(V,E)$, $0\leq|F|-|V|+k(F)\leq |E|-|V|+1$.
$(t^{-1}+1)^{|F|-|V|+k(F)}$ now can be expressed as
\begin{eqnarray*}
& &t^{-({|F|-|V|+k(F)})}+({|F|-|V|+k(F)})t^{-(|F|-|V|+k(F)-1)}+\\
& & \left(
\begin{array}{c}
|F|-|V|+k(F)\\
2
\end{array}
\right)t^{-(|F|-|V|+k(F)-2)}+\cdots.
\end{eqnarray*}

Note that $G$ is connected and bridgeless, we have
$|F|-|V|+k(F)=|E|-|V|+1$ if and only if $F=E$. Hence, we have
$n=-|E|+|V|-1$ and $a_n=(-1)^{|E|-|V|+1}$. Furthermore,
$|F|-|V|+k(F)=|E|-|V|$ if and only if $F=E-e$ for $e\in E$ or, by
Proposition \ref{111} $F=E-S$, where $S$ is a pairwise-disconnecting
set of $G$. Thus,

\begin{eqnarray*}
a_{n+1}&=&(-1)^{|E|-|V|+1}(|E|-|V|+1)+(-1)^{|E|-|V|}|E|+
\sum_{E-S}(-1)^{|E-S|-|V|+1}\\
&=&(-1)^{|E|-|V|+1}(-|V|+1)+\sum_{E-S_M}\sum_{S\subset S_M}(-1)^{|E-S|-|V|+1}\\
& & \textrm{(By Proposition \ref{222})}\\
&=&(-1)^{|E|-|V|+1}(-|V|+1+\sum_{S_M}(|S_M|-1)).
\end{eqnarray*}

$\Box$

\begin{remark}
Results of Lemmas \ref{hi} and \ref{lo} are dual in the sense that
(maximal) multiple edge corresponds to (maximal)
pairwise-disconnecting set by taking the dual when they are both
plane graphs.
\end{remark}

\begin{theorem}\label{12}
Let $G$ be a connected bridgeless and loopless positive graph. Then
the highest and lowest degrees of $Q[G]$ are $3|E|-2|V|+2$ and
$-|E|-2|V|+2$, respectively. Furthermore,
\begin{enumerate}
\item[(1)] the coefficient of the term with the highest degree is
$(-1)^{|E|-|V|+1}$,
\item[(2)] the coefficient of the term with the lowest degree is
$(-1)^{|V|-1}$,
\item[(3)] the coefficient of the term with the second-highest degree is
$(-1)^{|E|-|V|+1}(-|V|+1+\sum_{S_M}(|S_M|-1))$,
\item[(4)] the coefficient of the term with the second-lowest degree is $(-1)^{|V|-1}(|V|-1-|E|+\sum_{I_M}(|I_M|-1))$.
\end{enumerate}
\end{theorem}

\noindent{\bf Proof.} It follows from Lemmas \ref{this}, \ref{hi}
and \ref{lo}. $\Box$

\begin{remark}
Theorem \ref{12} (1) and (2) are the generalization of Lemma
\ref{degree} from planar graphs to all abstract (not necessarily
planar) graphs. It is not difficult to verify that when $G$ is a
plane graph, Theorem \ref{12} (1) and (2) coincide with Lemma
\ref{degree}.
\end{remark}

\section{The proof of Theorem \ref{eq}}

To prove Theorem \ref{eq}, we first need to further study the
properties of maximal pairwise-disconnecting sets. Let $G$ be a
connected bridgeless graph and $S$ be a pairwise-disconnecting set.
For any $e\in S$, we define $S_M(e)$ to be the union of $\{e\}$ and
the set of all bridges of $G-e$.

\begin{proposition}\label{333}
For any $e,f\in S$, $S_M(e)=S_M(f)$ and it is exactly the
unique maximal pairwise-disconnecting set containing $S$.
\end{proposition}

\noindent{Proof.} It suffices for us to prove that $S_M(e)\subset
S_M(f)$. It is clear that $e\in S_M(f)$ since $e,f\in S$ implying
that $\{e,f\}$ constitutes a 2-edge cut of $G$. For any $g\in
S_M(e)$ and $g\neq e$, $g$ is a bridge of $G-e$. Recall that
$\{e,f\}$ is a 2-edge cut of $G$ and suppose that $G-e-f=G_1\cup
G_2$ and $g\in E(G_2)$. $g$ is a bridge of $G-e$ implies that $g$ is
a bridge of $G_2$, and is also a bridge of $G-f$. Hence, $g\in
S_M(f)$ and we proved that $S_M(e)\subset S_M(f)$.

It is clear that $S\subset S_M(e)$. Now we prove that $S_M(e)$ is a
maximal pairwise-disconnecting set. According to the definition of
$S_M(e)$ we know that $k(G-S_M(e))=|S_M(e)|$. By Proposition
\ref{111}, we have $S_M(e)$ is a pairwise-disconnecting set. To
prove the maximality of $S_M(e)$, we suppose that $g\notin S_M(e)$
and $\{g\}\cup S_M(e)$ is a pairwise-disconnecting set. Then
$\{e,g\}$ is a 2-edge cut of $G$ and $g$ is a bridge of $G-e$,
contradicting $g\notin S_M(e)$. $\Box$

\begin{proposition}\label{444}
Any two distinct maximal pairwise-disconnecting sets of a connected
bridgeless graph are disjoint.
\end{proposition}

\noindent{Proof.} Suppose that $S_M^1$ and $S_M^2$ are two distinct
maximal pairwise-disconnecting sets of a connected bridgeless graph
$G$ and $e\in S_M^1\cap S_M^2$. By Proposition \ref{333}, $S_M^i$
($i=1,2$) will both be the union of $\{e\}$ and the set of all
bridges of $G-e$ and hence, will be equal, a contradiction. $\Box$

\begin{proposition}\label{555}
A pairwise-disconnecting set $S=\{e_1,e_2,\cdots,e_k\}$ of a
connected bridgeless graph $G$ as shown in Fig. \ref{cyc} is maximal
if and only if each $G_i$ ($i=1,2,\cdots,k$) is bridgeless.
\end{proposition}

\noindent{Proof.} It is obvious. $\Box$

Now we are in a position to prove Theorem \ref{eq}.

\vskip0.2cm

\noindent{\bf Proof.} If $\tilde{L}$ is not connected, then the $<$
of Eq. (\ref{e1}) holds. From the proof of Theorem \ref{lei}, we
know that $c(\tilde{L})=c(L)-2$ if and only if the Case 2a happens.

If $\tilde{L}$ is connected, then the $=$ of Eq. (\ref{e1}) holds.
Let $a_2$ (resp. $b_2$) be the coefficient of the degree $e(G)+2f(G)-4$ (resp. $-e(G)-2v(G)+8$) in
$Q[\tilde{G}]$. From
the proof of Theorem \ref{lei}, the equality of Theorem \ref{lei}
holds if and only if $a_2\neq 0$ and $b_2\neq 0$.  There are two
cases.

\noindent{\bf Case 1.} $G''=G/e$ is loopless.

Note that $G''$ is loopless means that $\mu(e)=1$. By Eq. (\ref{case1}) and Theorem \ref{12},
we obtain that

\begin{eqnarray*}
a_2&=&(-1)^{|E|-|V|+1}(-|V|+1+\sum_{S_M}(|S_M|-1))+(-1)^{|E''|-|V''|+1}-\\
& &(-1)^{|E''|-|V''|+1}(-|V''|+1+\sum_{S_M''}(|S_M''|-1))\\
&=&(-1)^{|E|-|V|+1}(-|V|+1+\sum_{S_M}(|S_M|-1)+1+|V''|-1-\sum_{S_M''}(|S_M''|-1))\\
&=&(-1)^{|E|-|V|+1}(\sum_{S_M}(|S_M|-1)-\sum_{S_M''}(|S_M''|-1)).
\end{eqnarray*}

Maximal pairwise-disconnecting sets of $G$ can be divided into two
classes: those containing the edge $e$ and those not containing the
edge $e$. Let $S_M$ be a maximal pairwise-disconnecting set of $G$.
By Proposition \ref{555}, we obtain that $G-S_M=G_1\cup
G_2\cup\cdots\cup G_{|S_M|}$ and each $G_i$ ($i=1,2,\cdots,|S_M|$)
is bridgeless. If $e\in S_M$, suppose that $e$ connects $G_i$ to
$G_{i+1}$ for some $i$. Since the one-point join of $G_i$ and
$G_{i+1}$ is bridgeless, by Proposition \ref{555} we have $S_M-e$ is
a maximal pairwise-disconnecting set of $G''$. If $e\notin S_M$,
suppose $e\in E(G_i)$ for some $i$. Since $G_i/e$ is bridgeless, we
have $S_M$ is also a maximal pairwise-disconnecting set of $G''$.

Conversely, let $S_M''$ be a maximal pairwise-disconnecting set of
$G''$. By Proposition \ref{555}, we obtain that $G''-S_M''=G_1''\cup
G_2''\cup\cdots\cup G_{|S_M''|}''$ and each $G_i''$
($i=1,2,\cdots,|S_M''|$) is bridgeless. Note that the two
end-vertices $u$ and $v$ of the edge $e$ of $G$ is identified to
become one vertex, say $u''$, in $G''$. Suppose that $u''\in
V(G_i'')$ for some $i$ and $G-S_M''=G_1''\cup \cdots\cup
G_{i-1}''\cup G_i\cup G_{i+1}''\cup\cdots\cup G_{|S_M''|}''$. Then
$G_i''=G_i/e$. If $e$ is not a bridge of $G_i$, then $S_M''$ is a
maximal pairwise-disconnecting set of $G$. If $e$ is a bridge of
$G_i$, then $S_M''\cup\{e\}$ will be a maximal
pairwise-disconnecting set of $G$. Suppose that $G$ has exactly $k$
maximal pairwise-disconnecting sets containing the edge $e$. Then
\begin{eqnarray*}
a_2&=&(-1)^{|E|-|V|+1}k.
\end{eqnarray*}

Furthermore, by Proposition \ref{444}, $k=0$ or $1$. Thus $a_2\neq
0$ iff $G$ has (a unique) maximal
pairwise-disconnecting set containing the edge $e$ iff $G-e$ has bridges.

Similarly, we have

\begin{eqnarray*}
b_2&=&(-1)^{|V|-1}(|V|-1-|E_s|)-(-1)^{|V''|-1}+(-1)^{|V''|-1}(|V''|-1-|E_s''|)\\
&=&(-1)^{|V|-1}(|V|-1-|E_s|+1-|V''|+1+|E_s''|)\\
&=&(-1)^{|V|-1}(|E_s''|+2-|E_s|).
\end{eqnarray*}

It is not difficult to see that $|E_s|=|E_s''|+1+|N(u)\cap N(v)|$.
So $b_2\neq 0$ iff $|N(u)\cap N(v)|\neq 1$.

Moreover, $G-e$ has bridges imply that $|N(u)\cap N(v)|\leq 1$ (see
Fig. \ref{cyc}). Thus $a_2\neq 0$ and $b_2\neq 0$ if and only if
$G-e$ has bridges and $N(u)\cap N(v)=\emptyset.$

\noindent{\bf Case 2.} $G''$ has loops.

This means $\mu(e)\geq 2$. For any edge $f\in E(G)$, which is
parallel to $e$, if $G-e-f$ is disconnected, then $\tilde{L}$ will
be a split link. Hence $G-e-f$ is connected. Recall that $\tilde{G}$
corresponding to $\tilde{L}$ is obtained from $G$ by changing the
sign of $e$ from $+$ and $-$. Note that $e$ and $f$ will cancel each
other in $\tilde{G}$ by the second Reidemeister move and
$\tilde{G}-e-f=G-e-f$ is positive and loopless, we have
$c(\tilde{L})=c(L)-2$ if and only if $G-e-f$ is connected and
bridgeless. $\Box$

\section{Examples and further discussions}

In this section, we first provide an example to illustrate Theorem
\ref{eq}. It is well known that rational knots are alternating and
by changing a crossing of a rational knot we still obtain a rational
knot.

\begin{figure}[htbp]
\centering
\includegraphics[width=3cm]{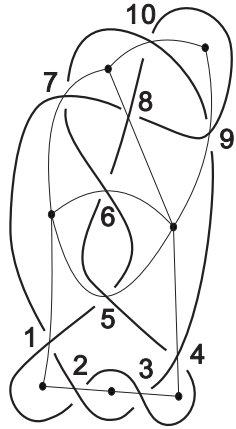}
\renewcommand{\figurename}{Fig.}
\caption{The knot $10_{14}$ (the dashed curve) and its corresponding
graph (the thick curve).}\label{exa}
\end{figure}

\begin{example}
The rational knot $10_{14}$ (see \cite{Ad} P. 47) (the dashed curve)
and its corresponding graph $G$ (the thick curve) are shown in Fig.
\ref{exa}. For $i=1,2,3,4$, $\mu(i)=1$, $G-i$ has bridges, the two
end-vertices of $i$ have no common neighbors. For $i=5,6$,
$\mu(i)=2$, $G-\{5,6\}$ is connected and bridgeless. For $i=7,8$,
$\mu(i)=1$, $G-i$ has no bridges. For $i=9,10$, $\mu(i)=1$, $G-i$
has bridges, the two end-vertices of $i$ have one common neighbor.
Hence the crossing number is reduced exactly by 2 after changing the
crossing $i$ for $i=1,2,3,4,5,6$ and reduced by 3 or more after
changing the crossing $i$ for $i=7,8,9,10$.

\end{example}

\vskip0.2cm

In the Dale Rolfsen's Knot table, if an alternating knot diagram
corresponds to a negative plane graph, we shall take its mirror
image to obtain a positive plane graph. Among alternating knots
whose crossing number is less than 10, there are only 11 knot
diagrams whose corresponding positive plane graphs satisfy
conditions of Corollary \ref{eq1}, and they are
$5_1$,$7_1$,$7_3$,$8_3$,$8_5$,$9_1$,$9_3$,$9_4$,$9_9$,$9_{10}$,$9_{35}$.
There are only 5 knot diagrams whose corresponding positive plane
graphs satisfy conditions of Corollary \ref{eq2}, and they are
$5_1$,$7_1$,$8_5$,$9_1$,$9_{35}$.

\begin{figure}[htbp]
\centering
\includegraphics[width=8.5cm]{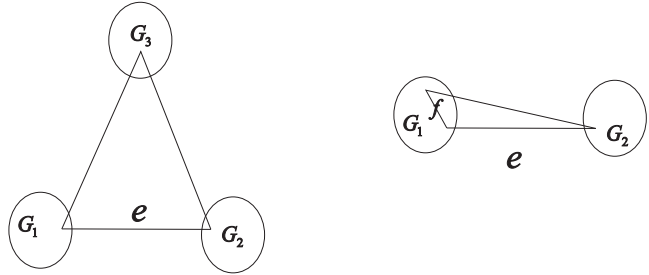}
\renewcommand{\figurename}{Fig.}
\caption{Two types of graphs with $\mu(e)=1$, $G-e$ has bridges and
$N(u)\cap N(v)\neq\emptyset$. In both types $G_i$ is bridgeless and
in the second type $f\in E(G_1)$.}\label{ty}
\end{figure}

Moreover, in graph theory, it is easy to judge whether an edge is a
bridge or not. As for the condition $N(u)\cap N(v)=\emptyset$, under
conditions $\mu(e)=1$ and $G-e$ has bridges, there are only two
types of graphs with $N(u)\cap N(v)\neq\emptyset$ as shown in Fig.
\ref{ty}.

Finally, although sufficient and necessary conditions of Theorem
\ref{eq} and two corollaries are very simple, applications of
Theorem \ref{eq} or its two corollaries are still very limited since
the properties, non-split and alternating of $\tilde{L}$, have not
been converted to conditions of $G$ (and the edge $e$). We pose the
following two problems for further study.

\noindent {\bf Problem 1.}  Let $G$ be a positive plane graph, $e$
be an edge of $G$. Let $\tilde{L}$ be the link whose diagram
corresponds to the plane graph obtained from $G$ by changing the
sign of $e$ from $+$ to $-$. We ask which conditions should be
satisfied by $G$ and $e$ to guarantee that the link $\tilde{L}$ is
non-split?

We note that {\bf Problem 1} appears in Page 143 of \cite{Ad} as an unsolved question.

\noindent {\bf Problem 2.}  Let $G$ be a positive plane graph, $e$
be an edge of $G$. Let $\tilde{L}$ be the link whose diagram
corresponds to the plane graph obtained from $G$ by changing the
sign of $e$ from $+$ to $-$. We ask which conditions should be
satisfied by $G$ and $e$ to guarantee that the link $\tilde{L}$ is
alternating?

\vskip0.2cm \noindent{\bf \large Acknowledgements} \vskip0.2cm

This paper was supported by NSFC Grant No. 10831001 and Grant No. 11271307.
We thank the referees for their suggestions.

\bibliographystyle{model1b-num-names}
\bibliography{<your-bib-database>}

\end{document}